\documentclass[12pt]{article}
\usepackage{amsfonts}

\newcommand{\sect}[1]{\setcounter{equation}{0}\section{#1}}
\newcommand{\subsect}[1]{\subsection{#1}}

\def\be{\begin{equation}}
\def\ee{\end{equation}}
\def\bea{\begin{eqnarray}}
\def\eea{\end{eqnarray}}

\def\1{\'{\i}}                           
\def\R{{\mathbb R}}
\def\Z{{\mathbb Z}}

\def\kk{K}
\def\pp{P}
\def\hh{H}
\def\cuno{C_1}
\def\cdos{C_2}
\def\dd{D}
\def\zzt{\tau}
\def\zzs{\sigma}
\def\tildezzt{\tau}

\def\casa{W_1}
\def\casb{W_2}

\def\casazzt{W_{1,\zzt}}
\def\casbzzt{W_{2,\zzt}}

\def\mma{\mu}
\def\mmc{\nu}

\begin{document}

\thispagestyle{empty}

\ 
\hfill\

\ 
\vspace{2cm}

\begin{center}
{\Large{\bf{New time-type and space-type 
 non-standard quantum   algebras  and   discrete symmetries}}} 
\end{center}

\bigskip\bigskip

\begin{center}   
 Francisco J. Herranz 
\end{center}

\begin{center} {\it { 
Departamento de F\1sica\\
Escuela Polit\'ecnica Superior\\ 
Universidad de Burgos\\
E-09006 Burgos, Spain }}
\end{center}

\bigskip\bigskip

\begin{abstract}
Starting from the classical $r$-matrix of the non-standard (or Jordanian)
quantum deformation of the $sl(2,\R)$ algebra,   new triangular quantum
deformations for the real Lie algebras $so(2,2)$, $so(3,1)$ and
$iso(2,1)$ are simultaneously constructed by using a graded contraction
scheme; these are realized as deformations of  conformal algebras of 
$(1+1)$-dimensional   spacetimes.  Time-type and space-type  quantum
algebras are considered according to the  generator that remains
primitive after deformation: either the time or the space translation,
respectively. Furthermore by introducing differential-difference
conformal realizations, these families of quantum   algebras are  shown
to be the symmetry algebras of either a  time or a space discretization
of $(1+1)$-dimensional (wave and Laplace)  equations on uniform lattices;
the relationship with  the known Lie symmetry approach to these discrete
equations is established by means of  twist maps.
\end{abstract}

\newpage


\sect{Introduction}

The non-standard (or Jordanian) quantum deformation of $sl(2,\R)\simeq
so(2,1)$ \cite{Ohn}, $U_z(sl(2,\R))$, has been the starting point in
the construction of  non-standard  (or triangular) quantum algebras in
higher dimensions. In particular, by taking two copies of $U_z(sl(2,\R))$
and applying the same procedure as in the standard (Drinfel'd--Jimbo)
case \cite{Ita}, a quantum $so(2,2)$ algebra has been obtained in
\cite{beyond}, while the corresponding deformation for $so(3,2)$ has been
found in \cite{vulpi}. These  quantum algebras have been realized as
deformations of conformal algebras for the Minkowskian spacetime.
Furthermore, by following either a contraction approach \cite{beyond} or
a deformation embedding method \cite{nulla}, non-standard quantum
deformations for other Lie algebras have been deduced; amongst them  it
is remarkable the appearance of a non-standard quantum Poincar\'e algebra,
which can be considered  as  a  quantum conformal  algebra for
the Carroll spacetime, or alternatively  and more interesting,  as a
null-plane quantum Poincar\'e algebra  \cite{nulla,nullb}.  All these
results are summarized in the following diagram where the vertical arrows
indicate the  contractions leading  to quantum Poincar\'e algebras:
$$
\begin{array}{ccccc}
U_z(sl(2,\R))&\longrightarrow & U_z(sl(2,\R))\oplus
U_{-z}(sl(2,\R))\simeq U_z(so(2,2))&\longrightarrow &
U_z(so(3,2))\\[5pt]
 \Big\downarrow\
\varepsilon
\to 0&&
\qquad\qquad\qquad\qquad\qquad\qquad\quad 
\Big\downarrow\ \varepsilon
\to 0& & \Big\downarrow\ \varepsilon \to 0\\[5pt]
U_z(iso(1,1))&\longrightarrow&
\quad\mbox{\footnotesize{Null-plane  Poincar\'e algebra}}\qquad
U_z(iso(2,1))&\longrightarrow & U_z(iso(3,1))\cr
 \end{array} 
$$

A first  aim of this paper is  to provide, starting again  from 
$U_z(sl(2,\R))$, a new way in the construction of  non-standard quantum
algebras obtaining a   new non-standard quantum $so(2,2)$ algebra which
could be the cornerstone of further constructions in higher dimensions. 
The essential idea is to require that  $U_z(sl(2,\R))$ remains as a Hopf
subalgebra, or to be more precise, to keep its underlying Jordanian 
classical  $r$-matrix, $r=z J_3\wedge J_+$, as the element generating the
whole  deformation for 
$so(2,2)$. Hence this approach can be seen as a kind of {\em complete}
deformation embedding method  leading to $U_z(sl(2,\R))\subset
U_z(so(2,2))$,  so that this seems to be  a more feasible and applicable
quantum deformation procedure than the involved  one used in \cite{vulpi}
for   $so(3,2)$ when the extension to higher dimensions is attacked.

Two    choices for such Jordanian  classical  $r$-matrix associated to
$so(2,2)$ naturally appear: one gives rise to  a  {\em time-type} quantum
deformation characterised by  a primitive generator of time translations,
meanwhile the other leads to   a {\em space-type} deformation determined
by a primitive generator of space translations; the Drinfel'd--Jimbo
counterpart of these types of deformations can be found in \cite{BHOS}.
Furthermore by using graded contractions this task is carried out for the
real Lie algebras  $so(2,2)$, $so(3,1)$  and $iso(2,1)$, simultaneously;
the quantum algebras so obtained are realized as deformations of 
conformal algebras of  $(1+1)$D   spacetimes.

The  second aim of this paper is to analyse the discrete symmetries
provided by both families of quantum algebras as differential-difference
conformal operators of either  a time  or a space discretization of some
$(1+1)$D differential equations (the wave and Laplace equations) on a
uniform lattice, and next to relate these results with the  Lie symmetry
analysis presented in \cite{Javier}. This objective is achieved by
following a similar procedure to the one used in \cite{schrod} with
respect to non-standard quantum Schr\"odinger algebras and their
associated discrete symmetries.

 The structure of the paper is as follows. We summarize in the next
section the  basic aspects of the $\Z_2\times \Z_2$ graded contractions
of  $so(2,2)$  in a conformal  basis as well  as their role of symmetry
algebras of  $(1+1)$D differential equations. The construction of the
time-type quantum deformation together with its  universal
 $R$-matrix is developed in the section 3. These quantum algebras are
shown to be the symmetry algebras of a time discretization of the wave
and Laplace equations on a uniform lattice in the section 4; the
relationship with the  Lie symmetry approach studied in \cite{Javier} is
 established by means of a  twist map. A parallel procedure with
 the space-type quatum deformation is carried out in the section 5.
Finally, an algebraic equivalence or {\em duality} between both types of
quantum algebras is introduced in the last section where we also comment
on their possible generalization to higher dimensions and the way of
obtaining new null-plane quantum Poincar\'e algebras.


\sect{Graded contractions of  $\mbox{\boldmath $so(2,2)$}$  and continuous
symmetries}

The $\Z_2\times \Z_2$  graded contractions of the real Lie algebra
$so(2,2)$ have been analysed in  \cite{beyond}, where  a distinguished
set of solutions has been explicitly considered and expressed in terms of
three contraction parameters $(\mu_1,\mu_2, \mu_3)$. Here we   shall
restrict ourselves to deal  with the most relevant contracted Lie
algebras  setting $(\mu_1,\mu_2, \mu_3)= (\mma,+1,\mmc)$,  so that all of
them are collectively    denoted $so_{\mma,\mmc}(2,2)$. Recall that each
contraction parameter can take either a positive, zero or  negative value
and whenever they are different from zero
 can  be scaled to  $\pm 1$.

At this dimension a generic Lie algebra in  the family
$so_{\mma,\mmc}(2,2)$ can be interpreted in two different frameworks:
either as the algebra of isometries of a $(2+1)$D spacetime (or a 3D
space), or as the  algebra of conformal transformations of a $(1+1)$D
spacetime (or a 2D space). In this paper we will adopt   the latter
interpretation, hence let us consider the   generators of  time
translations  $H$,  space translations $P$, boosts $K$, dilations $D$
and  special conformal transformations  $C_1$, $C_2$. In this basis,  the
Lie brackets of the set of graded contractions $so_{\mma,\mmc}(2,2)$ read
\be
\begin{array}{lll}
[K,H]=\mmc P&\qquad [K,P]=\mma H&\qquad [H,P]=0\cr
[D,H]=H&\qquad [D,C_1]=-C_1&\qquad [H,C_1]=-2\mmc D\cr
[D,P]=P&\qquad [D,C_2]=-C_2&\qquad [P,C_2]=2\mma D\cr
[K,C_1]=\mmc C_2&\qquad [K,C_2]=\mma C_1&\qquad [K,D]=0 \cr
[H,C_2]=2K&\qquad [P,C_1]=-2K&\qquad  [C_1,C_2]=0.
 \end{array} 
 \label{ba}
\ee
The two Casimirs of $so_{\mma,\mmc}(2,2)$ turn out to be
\be
\begin{array}{l}
\casa=K^2 +\mma\mmc D^2 -\frac 12  \mma(H C_1 + C_1 H) +\frac 12 \mmc (P
C_2 + C_2 P)\cr
\casb=K D + \frac 12 (H C_2 - C_1 P).
 \end{array} 
 \label{bb}
\ee

In what follows we  identify  each specific  real Lie algebra appearing
within the family $so_{\mma,\mmc}(2,2)$ (see table 1 below) and   comment
its physical (or  geometrical) role   according to  the (signs or zero)
values   of the pair
$(\mma,\mmc)$
\cite{beyond}: 
 
\noindent
$\bullet$  $so(2,2)$ when $(\mma,\mmc)\in\{(+,+), (-,-)\}$. This is the 
conformal algebra  of the $(1+1)$D Minkowskian spacetime;
alternatively, it can be seen as  the kinematical
algebra of the $(2+1)$D Anti-de Sitter spacetime.

\noindent
$\bullet$  $so(3,1)$ when $(\mma,\mmc)\in\{(+,-), (-,+)\}$. This is the
conformal algebra of the 2D Euclidean space  so that, under this
interpretation, $H$ should be considered as another generator of space
translations.  This algebra can also be  realized as the kinematical
algebra of the   $(2+1)$D  de Sitter spacetime.

\noindent
$\bullet$   $iso(2,1)$ when $(\mma,\mmc)\in\{(+,0), (0,+),(-,0),
(0,-)\}$. In the four cases, this is the kinematical algebra of the 
 $(2+1)$D  Minkowskian spacetime,  that is, the  $(2+1)$D Poincar\'e
algebra. However, this corresponds to the conformal algebra of the
$(1+1)$D Galilean spacetime whenever  $(\mma,\mmc)=(0,\pm)$, but to the
conformal algebra of the $(1+1)$D Carroll spacetime whenever 
$(\mma,\mmc)=(\pm,0)$.

\noindent
$\bullet$  $i'iso(1,1)$ when $(\mma,\mmc)=(0,0)$. This is the most
contracted algebra in the  family $so_{\mma,\mmc}(2,2)$ and has no known
conformal interpretation, although  is the algebra of
isometries of certain  3D space. Note that in this case $K$ is a central
generator.

The aforementioned conformal role  of the algebras $so_{\mma,\mmc}(2,2)$
(with the exception of $i'iso(1,1)$)  can be appreciated  more clearly 
by taking into account that: (i) The Lie brackets of the subalgebra
spanned by $\{K,H,P\}$  generate the  algebra of isometries of the
corresponding $(1+1)$D spacetime (or 2D Euclidean space). (ii) When the
dilation generator is added, we find the so called Weyl subalgebra
$\{K,H,P,D\}$ which is the similitude algebra of the 
$(1+1)$D spacetime. (iii) If  conformal transformations are also
considered, then we obtain the complete conformal Lie group  
$SO_{\mma,\mmc}(2,2)$; its quotient  with the subgroup generated by
$\{K,C_1,C_2,D\}$ is identified with the  $(1+1)$D conformal spacetime.

The relationship between $so_{\mma,\mmc}(2,2)$ and 
$(1+1)$D differential equations can be established by considering the 
usual (conformal)  vector field  representation   in terms of the space
and time coordinates $(x,t)$:
\be
\begin{array}{l}
 H=\partial_t\qquad P=\partial_x\qquad
K =- \mmc t \partial_x -\mma x\partial_t\qquad
D=-x\partial_x - t \partial_t\cr
C_1=(\mma x^2+\mmc t^2) \partial_t + 2\mmc x t   \partial_x\qquad
C_2=-(\mma x^2+\mmc t^2) \partial_x - 2 \mma x t   \partial_t  
 \end{array} 
 \label{bd}
\ee
where we exclude the degenerate case   $i'iso(1,1)$ with
$\mma=\mmc=0$. This is a zero-value  realization of the two Casimirs
(\ref{bb}).  The action of the Casimir of the  Lie subalgebra 
$\{K,H,P\}$, 
\be
E=\mmc P^2-\mma H^2  
\label{be}
\ee
 on a function $\Phi(x,t)$ through the 
representation (\ref{bd}) (choosing for $E$ the  zero eigenvalue) leads
to the following $(1+1)$D differential   equation:
\be
E\Phi(x,t)=0\quad
\Longrightarrow\quad\left(\mmc \frac{\partial^2}{\partial
x^2}-\mma \frac{\partial^2}{\partial t^2}
\right)\Phi(x,t)=0 .
\label{bf}
\ee
 We shall say that an operator ${\cal O}$ is a symmetry of the
equation $E\Phi(x,t)=0$ if ${\cal O}$ transforms solutions into
solutions, that is, $E{\cal O}=\Lambda E$ where $\Lambda$ is another
operator. Hence, $so_{\mma,\mmc}(2,2)$ is the symmetry algebra
of the  equation (\ref{bf}), since $E$ given by (\ref{be}) commutes with
$\{K,H,P\}$ and  in the realization (\ref{bd}) the remaining 
generators are also symmetry operators of (\ref{bf}) verifying
\be
 [E,D]=-2E\qquad
[E,C_1]= 4 \mmc t E\qquad [E,C_2]= - 4\mma  x E .
 \label{bg}
\ee
From this perspective, we find that the   equation (\ref{bf}) reproduces
the $(1+1)$D wave equation   when  the contraction parameters
$(\mma,\mmc)$ are either $(+,+)$ or $(-,-)$, which  in turn means that
$so(2,2)$ is its associated  algebra of  symmetry operators. Likewise,
$so(3,1)$ corresponding to  $(+,-)$ or $(-,+)$ arises as the symmetry
algebra of the 2D Laplace equation (in this case $t$ should be seen as
another space coordinate). Finally, the contraction with either $\mma=0$
or $\mmc=0$ leads to the 1D Laplace equation with   $iso(2,1)$ as its
symmetry Lie algebra.


\sect{Time-type quantum   algebras}

Let us   consider the   subalgebra of $so_{\mma,\mmc}(2,2)$ spanned
by  $\{D,H\}$ with Lie bracket $[D,H]=H$,  and the non-standard or
Jordanian classical
$r$-matrix given by \cite{Drinfelda,Drinfeldb}:\be
r=-\zzt\dd\wedge \hh
\label{ca}
\ee
which is a solution of the classical  Yang--Baxter equation and $\zzt$ is
the deformation parameter.
 As is well known the deformed commutator and coproduct
for this subalgebra     can be written as
\be
 [\dd,\hh]=\frac{{1- {\rm e}^{-\zzt\hh}}}{\zzt}\qquad
\Delta(\hh)=1\otimes \hh + \hh \otimes 1\qquad
\Delta(\dd)= 1\otimes \dd + \dd \otimes {\rm e}^{-\zzt\hh} .
\label{cb}
\ee
We recall that this structure is a Hopf subalgebra of
non-standard quantum  deformations of $sl(2,\R)$ 
\cite{Ohn,Demidov,Zakr,nonsb,Ogi},  $iso(1,1)$, $gl(2)$, $h_4$
and Schr\"odinger  algebras \cite{schrod};
this was also introduced in \cite{Majida,Majidb} in relation to an
approach to physics at the Planck scale.

 If we impose now 
the classical $r$-matrix (\ref{ca}) to be  the generating object of a
quantum deformation for the whole family $so_{\mma,\mmc}(2,2)$, then  the
cocommutator $\delta$ of a generator $X$ that defines the associated Lie
bialgebra is obtained  as 
$\delta(X)=[ 1\otimes X+X\otimes 1,r]$, namely,
\be
\begin{array}{l}
\delta(\hh)=0\qquad \delta(\dd)= -\zzt \dd\wedge\hh\cr
\delta(\pp)=\zzt \pp\wedge \hh\qquad
\delta(\kk)=-\zzt \mmc \dd\wedge \pp\cr
\delta(\cuno)=-\zzt \cuno\wedge \hh  \qquad
\delta(\cdos)=-\zzt \cdos\wedge \hh + 2 \zzt \dd \wedge
\kk .
\end{array} 
 \label{cc}
\ee
The coproduct $\Delta$ for the quantum 
  algebras denoted $U_\zzt(so_{\mma,\mmc}(2,2))$  is obtained by solving
the coassociativity condition $(1\otimes \Delta)\Delta=(\Delta\otimes
1)\Delta$, by requiring that (\ref{cb}) remains as a Hopf subalgebra of 
$U_\zzt(so_{\mma,\mmc}(2,2))$, and by taking into  account that $\delta$
is related to the first order of $\Delta$ on $\zzt$, $\Delta_{(1)}$, by 
$\delta=\Delta_{(1)}-\sigma\circ \Delta_{(1)}$ where 
$\sigma(X\otimes Y)=Y\otimes X$. The resulting coproduct turns
out to be
\be
\begin{array}{l}
\Delta(\hh)=1\otimes \hh + \hh \otimes 1\qquad
\Delta(\dd)= 1\otimes \dd + \dd \otimes {\rm e}^{-\zzt\hh}\cr
\Delta(\pp)= 1\otimes \pp + \pp \otimes {\rm e}^{\zzt\hh} 
\qquad\Delta(\cuno)= 1\otimes \cuno + \cuno \otimes
{\rm e}^{-\zzt\hh}\cr
\Delta(\kk)=1\otimes \kk + \kk\otimes 1 - 
\zzt \mmc \dd \otimes {\rm e}^{-\zzt\hh}\pp\cr
\Delta(\cdos)=1\otimes \cdos + \cdos
\otimes {\rm e}^{-\zzt\hh}  + 2 \zzt\dd \otimes {\rm e}^{-\zzt\hh}\kk
-\zzt^2 \mmc (\dd^2+\dd)\otimes {\rm e}^{-2\zzt\hh}\pp .
\end{array} 
 \label{cd}
\ee
Thereafter, the deformed commutation rules  are deduced by imposing
$\Delta$ to be an algebra homomorphism, that is,
$\Delta([X,Y])=[\Delta(X) ,\Delta(Y)]$; they are
\be
\begin{array}{l}
\displaystyle{[\kk,\hh]=\mmc {\rm e}^{-\zzt\hh}\pp\qquad 
[\kk,\pp]=\mma  \frac{ {\rm e}^{\zzt\hh}-1}{\zzt} 
\qquad [\hh,\pp]=0\qquad [\kk,\dd]=0}\cr
\displaystyle{  [\dd,\hh]=\frac{1- {\rm e}^{-\zzt\hh}}{\zzt}
\qquad [\dd,\cuno]=-\cuno +   \zzt \mmc \dd^2 \qquad
[\hh,\cuno]=-2\mmc \dd }\cr
 [\dd,\pp]= \pp\qquad [\dd,\cdos]=-\cdos \qquad
[\pp,\cdos]=2\mma \dd\cr
 [\kk,\cuno]=\mmc \cdos  \quad\ 
[\kk,\cdos]=\mma \cuno- \zzt \mma\mmc\dd^2\quad\ 
 [\hh,\cdos]={\rm e}^{-\zzt\hh}\kk +
\kk {\rm e}^{-\zzt\hh}\cr
[\pp,\cuno]=-2\kk - 
\zzt \mmc(\dd\pp+\pp\dd)\qquad [\cuno,\cdos]=-\zzt\mmc 
(\dd\cdos+\cdos\dd) .
\end{array} 
 \label{ce}
\ee
The deformed Casimirs of $U_\zzt(so_{\mma,\mmc}(2,2))$ are given by
\bea
&&\casazzt=K^2 +\mma\mmc D^2 -\frac 12 \mma\left( \frac{
 {\rm e}^{\zzt\hh}-1}{\zzt} \,   C_1 + C_1  \frac{ {\rm
e}^{\zzt\hh}-1}{\zzt}
\right) +\frac 12 \mmc (P C_2 + C_2 P)\cr
&&\qquad\qquad +\frac 12 \mma\mmc \left({\rm e}^{\zzt\hh}\dd^2 + \dd^2
{\rm e}^{\zzt\hh} \right)-\mma\mmc \dd^2\cr
&&\casbzzt=K D + \frac 12 \left(\frac {{\rm e}^{\zzt\hh}-1}{\zzt}\,  C_2
- C_1 P\right)+\frac 12 \zzt \mmc \dd^2 \pp.
 \label{cf}
\eea

Mathematical and physical properties of the  quantum algebras
$U_\zzt(so_{\mma,\mmc}(2,2))$ are characterized by their {\em primitive
generator} $H$ (that  with a vanishing cocommutator). More explicitly,
since the product $\zzt \hh$ has to be dimensionless in order to have a
homogeneous coproduct, the deformation parameter $\zzt$ has, in
principle,  the dimension of a {\em time} (notice that for
$U_\zzt(so(3,1))$, within the conformal interpretation,
$\zzt$ would be a length). This is similar  to  what happens with the
well known $\kappa$-Poincar\'e algebra 
\cite{Lukierskia,Giller,Lukierskib}, realized as a kinematical algebra
of the Minkowskian spacetime, and where the time translation generator is
also primitive; the deformation parameters $\kappa$ and $\zzt$ would be
related by $\kappa=1/\zzt$. Furthermore, as we shall show in the next
section, these time-type quantum algebras  directly lead   to a time
discretization of the   symmetries (\ref{bd}) and   equation (\ref{bf})
on a uniform  lattice.

On the other hand, at the level of Hopf subalgebras of
$U_\zzt(so_{\mma,\mmc}(2,2))$ two remarkable structures arise:

\noindent
$\bullet$  The generators
$\{\dd,\hh,\cuno\}$ close a Hopf subalgebra isomorphic either to the  
Jordanian  quantum $sl(2,\R)\simeq so(2,1)$ algebra if $\mmc\ne 0$, or
to a non-standard quantum Poincar\'e  algebra   $U_\zzt(iso(1,1))$  under
the contraction $\mmc=0$. 

\noindent
$\bullet$  The generators
$\{\kk,\hh,\pp,\dd\}$  span a Hopf subalgebra which is the similitude
algebra of a  $(1+1)$D spacetime. Therefore, as a byproduct of our
construction, we obtain for each member in the family
$U_\zzt(so_{\mma,\mmc}(2,2))$ a new quantum deformation of the Weyl
subalgebra of the corresponding conformal algebra. Notice that
$\{\kk,\hh,\pp\}$ only close a Hopf subalgebra whenever $\mmc=0$.

For the sake of clarity the specific Hopf subalgebras 
spanned by $\{\dd,\hh,\cuno\}$ and $\{\kk,\hh,\pp,\dd\}$ for
each quantum algebra in
the family $U_\zzt(so_{\mma,\mmc}(2,2))$ are displayed in the  table 1
according to the values of the pair $(\mma,\mmc)$. The horizontal
arrows indicate the contraction $\mma=0$  and the vertical ones the
contraction $\mmc=0$.  The symbols  ${\cal
{WM}}$, ${\cal {WE}}$, ${\cal {WG}}$ and ${\cal {WC}}$ mean, in this
order, the Weyl subalgebra of the Minkowskian, Euclidean, Galilean and
Carroll planes,  thus  
  reminding   the corresponding conformal spaces, while ${\cal
{WA}}$ means the Abelian algebra enlarged with a dilation generator. 
In this context, we
remark that other non-standard quantum deformations of these Weyl
algebras have been carried out in \cite{weyl}, the   underlying classical
$r$-matrix of which reads in our notation  $r= \omega(K\wedge H + D\wedge
P)$;  their generalization to  higher
dimensions can be found in \cite{vulpi}.
 We also  recall
that other non-standard classical $r$-matrices for $so(3,2)$ and
$so(4,2)$ (expressed as conformal algebras) can be
found in \cite{Lukierskic}.

\bigskip

{\small
\noindent
{\bf {Table 1}}. Hopf subalgebras $\{\dd,\hh,\cuno\}$ and
$\{\kk,\hh,\pp,\dd\}$, and associated  (difference and/or differential)
equations of the time-type quantum algebras
$U_\zzt(so_{\mma,\mmc}(2,2))$.
$$
\begin{array}{ccccc}
\hline\\[-8pt]
(+,+)\quad U_\zzt(so(2,2))&\longrightarrow &(0,+)\quad U_\zzt(iso(2,1))
&\longleftarrow &(-,+)\quad U_\zzt(so(3,1))\\[5pt]
U_\zzt(sl(2,\R)) \quad U_\zzt({\cal {WM}})& &
U_\zzt(sl(2,\R)) \quad U_\zzt({\cal {WG}})
& &U_\zzt(sl(2,\R)) \quad  U_\zzt({\cal {WE}})\\[5pt]
(\partial^2_x- \Delta^2_t)\Phi =0
& &\partial^2_x \Phi =0
& &(\partial^2_x+ \Delta^2_t)\Phi =0\\[5pt]
\downarrow& &\downarrow & &\downarrow\\[5pt]
(+,0)\quad U_\zzt(iso(2,1))&\longrightarrow &(0,0)\quad
U_\zzt(i'iso(1,1)) &\longleftarrow &(-,0)\quad U_\zzt(iso(2,1))\\[5pt]
U_\zzt(iso(1,1)) \quad U_\zzt({\cal {WC}})
& &U_\zzt(iso(1,1)) \quad U_\zzt({\cal {WA}})& &
U_\zzt(iso(1,1))\quad  U_\zzt({\cal {WC}})\\[5pt]
 \Delta^2_t\Phi =0
& &\mbox{Degenerate equation} 
& & \Delta^2_t\Phi =0\\[5pt]
\uparrow& &\uparrow & &\uparrow\\[5pt]
(+,-)\quad U_\zzt(so(3,1))&\longrightarrow &(0,-)\quad U_\zzt(iso(2,1))
&\longleftarrow &(-,-)\quad U_\zzt(so(2,2))\\[5pt]
U_\zzt(sl(2,\R)) \quad U_\zzt({\cal {WE}})& &
U_\zzt(sl(2,\R)) \quad U_\zzt({\cal {WG}})
 & &U_\zzt(sl(2,\R)) \quad  U_\zzt({\cal {WM}})\\[5pt]
(\partial^2_x+ \Delta^2_t)\Phi =0
& &\partial^2_x \Phi =0
& &(\partial^2_x- \Delta^2_t)\Phi =0\\[5pt]
\hline 
 \end{array}
$$
}


\subsect{Universal quantum  $\mbox{\boldmath $R$}$-matrix}

Different constructions of the  universal quantum $R$-matrix associated
to the non-standard quantum  deformation of the Borel algebra (of the
type  $[D,H]=H$) have  appeared in the
literature \cite{nonsb,Ogi,Shariati,BH}, mainly in relation to the
Jordanian  quantum $sl(2,\R)$ algebra.   If we consider the quantum Borel
algebra written in the  form of  (\ref{cb}), then the  universal 
$R$-matrix turns out to be
\cite{BH}
\be
{\cal R}=\exp\left\{\tau H\otimes D\right\}
\exp\left\{-\tau D\otimes H\right\}  
\label{cg}
\ee
 which  is a solution of the quantum Yang--Baxter equation and also
fulfils
\be
{\cal R}\Delta(X){\cal R}^{-1}=\sigma\circ\Delta(X)\quad\mbox{for}\quad X\in
\{\hh,\dd\} .
\label{ch}
\ee
  As it could be expected, the element
 (\ref{cg}) is  a  triangular universal  $R$-matrix  for the {\em whole}
family $U_\zzt(so_{\mma,\mmc}(2,2))$ since it contains (\ref{cb}) as a
Hopf subalgebra and the four remaining generators also verify (\ref{ch}):
\bea
&&\exp\left\{-\tau D\otimes H\right\}\Delta(\cuno)\exp\left\{\tau
D\otimes H\right\}=1\otimes \cuno +\cuno\otimes 1 + 2\zzt\mmc
\dd\otimes\dd\equiv f\cr &&\exp\left\{\tau H\otimes D\right\} f
\exp\left\{-\tau H\otimes D\right\}=\sigma\circ\Delta(\cuno) 
\label{ci}
\eea
\be
\begin{array}{l}
\exp\left\{-\tau D\otimes H\right\}\Delta(X)\exp\left\{\tau D\otimes
H\right\}=1\otimes X +X\otimes 1  \equiv \Delta_0(X)\cr
 \exp\left\{\tau H\otimes D\right\}  \Delta_0(X) \exp\left\{-\tau H\otimes
D\right\}=\sigma\circ\Delta(X)\quad \mbox{for}\quad X\in\{\pp,\kk,\cdos\}.
\end{array} 
\label{cj}
\ee

The lower dimensional  matrix  representation of
$U_\zzt(so_{\mma,\mmc}(2,2))$ is given by the following $4\times 4$ real
matrices:
\bea
&&H= \left(\begin{array}{cccc}
\frac 12\, \tildezzt\mmc&-\frac 12\, \tildezzt\mmc&-\mmc &0\\[2pt]
\frac 12\, \tildezzt\mmc&-\frac 12\,  \tildezzt\mmc&-\mmc&0\\[2pt]
1&-1&0&0\cr 0&0&0&0\end{array}\right) \qquad
P= \left(\begin{array}{cccc}
0&0&0 &\mma\cr 0&0&0&\mma\cr 0&0&0&0\cr 
1&-1&0&0\end{array}\right) \nonumber\\[5pt] 
&&K=\left(\begin{array}{cccc}
0&0&0 &0\cr 0&0&0&0\cr 0&0&0&\mma\cr  0&0&\mmc&0\end{array}\right) \qquad
D=\left(\begin{array}{cccc}
0&1&0 &0\cr 1&0&0&0\cr 0&0&0&0\cr  0&0&0&0\end{array}\right)\label{ck}
\\[5pt] 
&&C_1=  \left(\begin{array}{cccc}
\tildezzt\mmc&0&-\mmc &0\cr 0&\tildezzt\mmc&\mmc&0\cr 1&1&0&0\cr 
0&0&0&0\end{array}\right)\qquad C_2=  \left(\begin{array}{cccc}
0&0&0 &\mma\cr 0&0&0&-\mma\cr 
0&0&0&0 \cr
1&1&0&0\end{array}\right) .\nonumber
\eea
 We exclude the most contracted quantum algebra with $\mma=\mmc=0$,
since in this case the matrix of $K$ is degenerate (and this is a
central generator). This representation allows us to deduce a  $16\times
16$ matrix expression for  ${\cal R}$. Let us denote
 $\mbox{\boldmath $1$}$ and $\mbox{\boldmath $0$}$ the $4\times 4$ unit
and zero matrices;   under the representation (\ref{ck}) we find that 
$H^3={\mbox{\boldmath $0$}}$ so that  the quantum ${\cal R}$-matrix 
(\ref{cg})  reduces to
\be
\begin{array}{l}
 {\cal R}=( \mbox{\boldmath $1$}\otimes \mbox{\boldmath $1$}+\zzt
H\otimes D +\frac 12\zzt^2 H^2\otimes D^2)
( \mbox{\boldmath $1$}\otimes \mbox{\boldmath $1$}-\zzt D\otimes H
+\frac 12\zzt^2 D^2\otimes H^2) 
\end{array}
\label{cl}
\ee
which can finally be  written in block-matrix form as ${\cal R}=$
\be
{\small
\left(
\begin{array}{cccc|cccc|cccc|c}
1-\tildezzt^2\mmc&\tildezzt^2\mmc&0&0& 0&0&\tildezzt\mmc&0&
0&-\tildezzt\mmc&0&0& \cr 
0&1&0&0& -\tildezzt^2\mmc&\tildezzt^2\mmc&\tildezzt\mmc&0&
-\tildezzt\mmc&0&0&0&{\mbox{\boldmath $0$}}\cr 
0&0&1&0& -\tildezzt&\tildezzt&0&0& 0&0&0&0& \cr
0&0&0&1& 0&0&0&0& 0&0&0&0& \cr
\cline{1-13}
-\tildezzt^2\mmc&\tildezzt^2\mmc&\tildezzt\mmc&0& 1&0&0&0&
0&-\tildezzt\mmc&0&0& \cr
 0&0&\tildezzt\mmc&0& -\tildezzt^2\mmc&1+\tildezzt^2\mmc&0&0&
-\tildezzt\mmc&0&0&0&{\mbox{\boldmath $0$}}\cr 
-\tildezzt&\tildezzt&0&0& 0&0&1&0& 0&0&0&0& \cr
0&0&0&0& 0&0&0&1& 0&0&0&0& \cr
\cline{1-13}
 0&\tildezzt&-\tildezzt^2\mmc&0& 0&-\tildezzt&\tildezzt^2\mmc&0& 1&0&0&0&
\cr
\tildezzt&0&-\tildezzt^2\mmc&0& -\tildezzt&0&\tildezzt^2\mmc&0&
0&1&0&0&{\mbox{\boldmath $0$}}\cr 
0&0&0&0& 0&0&0&0& 0&0&1&0& \cr
0&0&0&0& 0&0&0&0& 0&0&0&1& \cr
\cline{1-13}
 \multicolumn{4}{c|}{\mbox{\boldmath $0$}}&
\multicolumn{4}{c|}{\mbox{\boldmath
$0$}}& \multicolumn{4}{c|}{\mbox{\boldmath $0$}}&{\mbox{\boldmath $1$}}\cr
\end{array}\right)
}
\label{cm}
\ee

So far  these last results could be further exploited in different
directions.  By one hand, the expression  (\ref{cg}) should allow one 
to obtain  a triangular $R$-matrix solution of the coloured Yang--Baxter
equation (that is, with spectral parameters); indeed this is 
formally rather similar to the universal $R$-matrix of the Jordanian
$gl(2)$ algebra used in \cite{Preeti} in order to deduce its coloured
realization. On the other hand, the matrices  (\ref{ck}) and   (\ref{cm})
could be applied in the computation of  the differential calculus on the 
quantum conformal spaces associated to $U_\zzt(so_{\mma,\mmc}(2,2))$. In
this respect see, for instance, \cite{Chang} where the construction of
the quantum Anti-de Sitter space  from the quantum algebra   
$SO_q(3,2)$ (of Drinfel'd--Jimbo type) has been carried out. The explicit
presence of the  contraction parameters  would enable a simultaneous
study of these problems for $so(2,2)$, $so(3,1)$ and
$iso(2,1)$ with a built-in scheme of contractions.


\sect{Discrete time symmetries}

Non-standard quantum Schr\"odinger algebras have been recently shown to
be the Hopf  algebras of  symmetries of a time (or a space)
discretization of the heat-Schr\"odinger equation on a uniform lattice
\cite{schrod}; in that construction the deformation parameter plays the
role of the  time (or space) lattice constant. Furthermore, by making use
of twist maps those discrete Schr\"odinger symmetries obtained from
quantum  algebras have been connected with the discretization (in a  
single variable) of the heat-Schr\"odinger equation deduced in
\cite{Luismi} by following  the usual Lie symmetry theory. In this
context, the remarkable point  is that   the same classical procedure has
been also applied  in \cite{Javier} to the study of the  symmetries of a
discretization of the  $(1+1)$D wave equation in both coordinates $(x,t)$
on a uniform lattice, showing that they are  difference   operators  
preserving the Lie algebra $so(2,2)$ as in the continuous case. Therefore
some kind of connection between the  results of \cite{Javier} and the 
quantum  $so(2,2)$ algebra here presented should exist  as it was already
established for discrete Shr\"odinger equations and quantum algebras in
\cite{schrod}.

Henceforth we follow a
parallel procedure  with the family $U_\zzt(so_{\mma,\mmc}(2,2))$ in
two steps. We first  introduce a differential-difference
realization showing that indeed this provides   discrete symmetries  of
a  time discretization of the equation (\ref{bf}). Secondly  we  give a
twist map that turns the deformed commutation rules of
$U_\zzt(so_{\mma,\mmc}(2,2))$ into the Lie commutators of
$so_{\mma,\mmc}(2,2)$ but keeping a (deformed) non-cocommutative
coproduct in such a manner that   a direct relationship with the  time
discretization of the wave equation  studied in \cite{Javier} from the
Lie symmetry approach can finally be   established.

 A differential-difference realization of   $U_\zzt(so_{\mma,\mmc}(2,2))$, 
which under the limit $\tau\to 0$ gives the continuous conformal realization
(\ref{bd}), reads
\bea
&&\hh=\partial_t \qquad \pp=\partial_x\nonumber\\[2pt]
&&\kk= -\mmc t  {\rm e}^{-\zzt\partial_t}  \partial_x - \mma x\left(
\frac{ {\rm e}^{\zzt\partial_t}-1 }{\zzt}\right) \qquad
\dd =- x \partial_x  
- t \left( \frac{1-{\rm e}^{-\zzt\partial_t}}{\zzt}\right)   \cr 
&&\cuno=(\mma x^2 +\mmc t^2 {\rm e}^{- \zzt\partial_t})\left( \frac{
{\rm e}^{\zzt\partial_t}-1}{\zzt}\right)
+ 2 \mmc x t \partial_x +\zzt \mmc (x \partial_x  
+  x^2 \partial^2_{x})\cr
&&\cdos=-(\mma x^2+\mmc t^2{\rm e}^{-2\zzt\partial_t})\partial_x
 - 2\mma x t \left( \frac{1- {\rm e}^{-\zzt\partial_t} }{\zzt} \right)
+\zzt \mmc t  {\rm e}^{-2 \zzt\partial_t}\partial_x .
\label{da}
\eea
In terms of (\ref{da}) both  deformed Casimirs (\ref{cf}) vanish.
The generators  $\{\kk,\hh,\pp\}$ close  a deformed  
subalgebra,  the Casimir of which is given by
\be
E_\zzt=\mmc \pp^2 -\mma \left( \frac{{\rm e}^{\zzt\hh}-1}{\zzt}\right)^2.
\label{db}
\ee
By introducing the realization (\ref{da})  we find a time discretization
of the  equation (\ref{bf})  on a uniform lattice with $x$ as a
continuous variable:
\be
E_\zzt\Phi(x,t)=0\quad
\Longrightarrow\quad \left\{
\mmc\frac{\partial^2}{\partial x^2}-
\mma \left(\frac{{\rm e}^{\zzt\partial_t}-1}{\zzt}\right)^2 
\right\}\Phi(x,t)=0 .
\label{dc}
\ee
The generators (\ref{da}) are symmetry
operators of (\ref{dc})   fulfilling
\bea
&& [E_\zzt,X]=0\quad \mbox{for}\quad X\in\{\kk,\hh,\pp\}
\qquad
[E_\zzt,\dd]=-2E_\zzt \cr
&&[E_\zzt,\cuno]= 4\mmc (t +\zzt + \zzt x \partial_x) E_\zzt
\qquad [E_\zzt,\cdos]= - 4\mma x E_\zzt .
 \label{dd}
\eea
 Hence  we conclude that $U_\zzt(so_{\mma,\mmc}(2,2))$ is the symmetry
algebra of the discrete  equation (\ref{dc}).

Next, let us consider the  so called minimal twist map, first introduced
in \cite{Abde} for the Jordanian quantum $sl(2,\R)$ algebra (here with
generators $\{\dd,\hh,\cuno\}$ and $\mmc\ne 0$) and  also used in
\cite{schrod} with other non-standard quantum algebras. This map  can be
implemented in  the whole family $U_\zzt(so_{\mma,\mmc}(2,2))$  as
\bea
&&{\cal H}=\frac{{\rm e}^{\zzt H}-1}{\zzt}\qquad {\cal P}=P\qquad {\cal
K}=K\qquad {\cal D}=D\cr &&{\cal C}_1=C_1-\zzt \mmc D^2\qquad {\cal
C}_2=C_2 .
\label{de}
\eea
These new generators verify the classical commutation rules (\ref{ba}),
while the coproduct remains deformed as
\be
\begin{array}{l}
 \displaystyle{\Delta({\cal H})=1\otimes {\cal H}+{\cal H}\otimes 1 +\zzt
{\cal H}\otimes {\cal H}}\qquad
 \displaystyle{\Delta({\cal P})=1\otimes {\cal P}+{\cal P}\otimes 1 +\zzt
{\cal P}\otimes {\cal H}}\\[5pt]
 \displaystyle{\Delta({\cal D})=1\otimes {\cal D}+{\cal D}\otimes
\frac{1}{1+\zzt {\cal H}} }\qquad
 \displaystyle{\Delta({\cal K})=1\otimes {\cal K}+{\cal K}\otimes 1-\zzt
\mmc {\cal D}\otimes \frac{{\cal P}}{1+\zzt {\cal H}}}\\[8pt]
 \displaystyle{\Delta({\cal C}_1)=1\otimes {\cal C}_1+{\cal C}_1\otimes
 \frac{1}{1+\zzt {\cal H}}-2\zzt \mmc {\cal D}\otimes \frac{1}{1+\zzt
{\cal H}}\,{\cal D}+ \zzt \mmc ({\cal D}^2+{\cal D})\otimes \frac{\zzt
{\cal H}}{(1+\zzt {\cal H})^2} }\\[8pt]
\displaystyle{ \Delta({\cal C}_2)=1\otimes {\cal C}_2+{\cal C}_2\otimes
\frac{1}{1+\zzt {\cal H}}+2\zzt {\cal D}\otimes \frac{1}{1+\zzt {\cal
H}}\,{\cal K}-
 \zzt^2 \mmc ({\cal D}^2+{\cal D})\otimes \frac{{\cal P}}{(1+\zzt {\cal
H})^2} }.
\end{array}
\label{df}
\ee
 The new generator   ${\cal H}$ verifies $\Delta((1 +\zzt {\cal H})^a)=(1
+\zzt {\cal H})^a\otimes (1 +\zzt {\cal H})^a$ for any real number $a$,
since $(1 +\zzt {\cal H})={\rm e}^{\zzt H}$.  

We apply the twist map (\ref{de}) to the   realization
(\ref{da}) and introduce the time shift operator $T_t={\rm e}^{\zzt
\partial_t}$ and the time difference operator $\Delta_t=(T_t-1)/\zzt$,
thus   finding  
\be
\begin{array}{l}
{\cal H}=\Delta_t\qquad {\cal P}=\partial_x\\[5pt]
 {\cal K}= -\mmc t  T_t^{-1}  \partial_x - \mma x \Delta_t\qquad
{\cal D} =- x \partial_x  
- t T_t^{-1} \Delta_t  \\[5pt]
 {\cal C}_1 =(\mma x^2 +\mmc t^2 T_t^{-2} )\Delta_t
 + 2 \mmc x t T_t^{-1} \partial_x  -   \zzt \mmc t T_t^{-2}\Delta_t  
\\[5pt]
 {\cal C}_2=-(\mma x^2+\mmc t^2 T_t^{-2} )\partial_x
 - 2\mma x t T_t^{-1}\Delta_t 
+\zzt \mmc t T_t^{-2} \partial_x .
\end{array}
\label{dg}
\ee
In this new basis  the Casimir of the   subalgebra $\{{\cal
K},{\cal H},{\cal P}\}$ is the undeformed one (\ref{be})
\be
{\cal E}=\mmc {\cal P}^2-\mma {\cal H}^2
\label{ddgg}
\ee
that written through (\ref{dg}) leads again to the
discrete equation (\ref{dc}):
\be
( \mmc \partial^2_x-
\mma\Delta^2_t)\Phi(x,t)=0 .
\label{dh}
\ee
The  new generators (\ref{dg}) are symmetry
operators of (\ref{dh})  now verifying   
\bea
 && [{\cal E} ,X]=0\quad \mbox{for}\quad X\in\{{\cal K},{\cal H},{\cal
P}\}
\qquad
[{\cal E} ,{\cal D}]=-2{\cal E}  \cr
&&[{\cal E},{\cal C}_1]= 4\mmc  t  T_t^{-1} {\cal E} 
\qquad [{\cal E} ,{\cal C}_2]= - 4\mma x {\cal E}  .
 \label{di}
\eea
In this way the relationship between $U_\zzt(so_{\mma,\mmc}(2,2))$ and 
 the symmetries of a time discretization of the wave equation deduced
from the Lie theory in \cite{Javier} clearly arises. In particular, let
us denote our generators, contraction parameters and variables by 
\bea 
&&\{{\cal H},{\cal P}, {\cal K}, {\cal D},{\cal
C}_1, {\cal C}_2\}= \{P_k,P_n,-L,-D,C_k,C_n\}\cr
&& (\mma,\mmc)=(s^2,+1)\qquad x= n\sigma \qquad   t=k\tau 
\label{dj}
\eea
with $s\ne 0$. If we  perform  the limits $n\to \infty$  and $\sigma\to
0$, subjected to the condition  $n\sigma=x$,  in the results given in
\cite{Javier} for $m=0$ (this implies that
$\Delta_n\to \partial_x$, $T_n\to 1$, $\Delta_k=\Delta_t$, $T_k=T_t$),
 then we recover the realization (\ref{dg}) and the discrete equation
(\ref{dh}). Consequently,  $U_\zzt(so_{s^2,1}(2,2))$ is the quantum
symmetry  algebra  of such equation and the deformation parameter $\tau$ 
is identified with the time lattice constant in the $t$ coordinate; the
space $x$ remains as a continuous variable. Recall that the solutions of
(\ref{dh}) has  also been obtained in
\cite{Javier}.

We  write down in the table 1 the particular equation (\ref{dh}) that
appears for each quantum algebra in the
 family $U_\zzt(so_{\mma,\mmc}(2,2))$.  It is worth noting that this
collective  treatment enables a clear view of the contraction limits
between  these (difference and/or differential) equations together with
their associated symmetry algebras.


\sect{Space-type quantum   algebras and\\
discrete space symmetries}

A second natural choice for a non-standard classical $r$-matrix 
for $so_{\mma,\mmc}(2,2)$, instead of (\ref{ca}), is to take
 \be
r=-\zzs\dd\wedge \pp
\label{ea}
\ee 
 where $\zzs$ is now the deformation parameter. 
If we follow the same steps described in section 3, we obtain a family of
 quantum algebras $U_\zzs(so_{\mma,\mmc}(2,2))$ characterised by a
primitive generator $P$ (instead of $H$). The resulting coproduct,
commutation rules and universal quantum $R$-matrix  are given by
\be
\begin{array}{l}
\Delta(\pp)=1\otimes \pp + \pp \otimes 1\qquad
\Delta(\dd)= 1\otimes \dd + \dd \otimes {\rm e}^{-\zzs\pp}\cr
\Delta(\hh)= 1\otimes \hh + \hh \otimes {\rm e}^{\zzs\pp}
\qquad \Delta(\cdos)=1\otimes \cdos + \cdos
\otimes {\rm e}^{-\zzs\pp}  \cr
\Delta(\kk)=1\otimes \kk + \kk\otimes 1 - 
\zzs \mma\dd \otimes {\rm e}^{-\zzs\pp}\hh\cr
\Delta(\cuno)= 1\otimes \cuno + \cuno \otimes {\rm e}^{-\zzs\pp} 
- 2 \zzs\dd \otimes {\rm e}^{-\zzs\pp}\kk +\zzs^2\mma (\dd^2+\dd)\otimes
{\rm e}^{-2\zzs\pp}\hh\cr 
\end{array}
\label{eb}
\ee
\be
\begin{array}{l}
\displaystyle{[\kk,\hh]=\mmc\,\frac{ {\rm e}^{\zzs\pp}-1}{\zzs}\qquad 
[\kk,\pp]=\mma {\rm e}^{-\zzs\pp}\hh
\qquad [\hh,\pp]=0\qquad [\kk,\dd]=0}\cr
\displaystyle{[\dd,\hh]=\hh\qquad 
[\dd,\cuno]=-\cuno\qquad
[\hh,\cuno]=-2\mmc\dd}\cr 
\displaystyle{ [\dd,\pp]= \frac{1- {\rm e}^{-\zzs\pp}}{\zzs}\qquad 
[\dd,\cdos]=-\cdos - \zzs \mma\dd^2\qquad
[\pp,\cdos]=2\mma\dd}\cr
\displaystyle{[\kk,\cuno]=\mmc \cdos +\zzs \mma\mmc\dd^2 \quad 
[\kk,\cdos]=\mma\cuno\quad  
[\pp,\cuno]=-{\rm e}^{-\zzs\pp}\kk - \kk
{\rm e}^{-\zzs\pp}}\cr
\displaystyle{ 
[\hh,\cdos]=2\kk + 
 \zzs\mma(\dd\hh+\hh\dd) \qquad  [\cuno,\cdos]=-\zzs\mma
(\dd\cuno+\cuno\dd)}  
\end{array}
\label{ec}
\ee
\be
{\cal R}=\exp\left\{\zzs P\otimes D\right\}
\exp\left\{-\zzs D\otimes P\right\}  .
\label{eec}
\ee
 At the level of Hopf subalgebras of  $U_\zzs(so_{\mma,\mmc}(2,2))$, we
find that  the generators $\{\dd,\pp,\cdos\}$ give rise  to  either a  
quantum $sl(2,\R)$ algebra if $\mma\ne 0$ or to a  quantum
  $iso(1,1)$  algebra    if $\mma=0$, meanwhile     $\{\kk,\hh,\pp,\dd\}$
close again a quantum Weyl  algebra; these  Hopf subalgebras are 
indicated in the table 2 for each pair $(\mma,\mmc)$.

\bigskip

{\small
\noindent
{\bf {Table 2}}. Hopf subalgebras $\{\dd,\pp,\cdos\}$ and
 $\{\kk,\hh,\pp,\dd\}$, and associated (difference and/or differential)
equations of the space-type quantum algebras
$U_\zzs(so_{\mma,\mmc}(2,2))$.
$$
\begin{array}{ccccc}
\hline\\[-8pt]
(+,+)\quad U_\zzs(so(2,2))&\longrightarrow &(0,+)\quad U_\zzs(iso(2,1))
&\longleftarrow &(-,+)\quad U_\zzs(so(3,1))\\[5pt]
U_\zzs(sl(2,\R)) \quad U_\zzs({\cal {WM}})& &
  U_\zzs(iso(1,1)) \quad U_\zzs({\cal {WG}})
& &U_\zzs(sl(2,\R)) \quad  U_\zzs({\cal {WE}})\\[5pt]
(\Delta^2_x-\partial^2_t)\Phi =0
& &\Delta^2_x \Phi =0
& &(\Delta^2_x+\partial^2_t)\Phi =0\\[5pt]
\downarrow& &\downarrow & &\downarrow\\[5pt]
 (+,0)\quad U_\zzs(iso(2,1))&\longrightarrow &(0,0)\quad
U_\zzs(i'iso(1,1)) &\longleftarrow &(-,0)\quad U_\zzs(iso(2,1))\\[5pt]
U_\zzs(sl(2,\R))\quad U_\zzs({\cal {WC}})
& &U_\zzs(iso(1,1)) \quad U_\zzs({\cal {WA}})& &
 U_\zzs(sl(2,\R)) \quad  U_\zzs({\cal {WC}})\\[5pt]
 \partial^2_t\Phi =0
& &\mbox{Degenerate equation} 
& & \partial^2_t\Phi =0\\[5pt]
\uparrow& &\uparrow & &\uparrow\\[5pt]
(+,-)\quad U_\zzs(so(3,1))&\longrightarrow &(0,-)\quad U_\zzs(iso(2,1))
&\longleftarrow &(-,-)\quad U_\zzs(so(2,2))\\[5pt]
U_\zzs(sl(2,\R)) \quad U_\zzs({\cal {WE}})& &
U_\zzs(iso(1,1)) \quad U_\zzs({\cal {WG}})
 & &U_\zzs(sl(2,\R)) \quad  U_\zzs({\cal {WM}})\\[5pt]
(\Delta^2_x+\partial^2_t)\Phi =0
& &\Delta^2_x \Phi =0
& &(\Delta^2_x-\partial^2_t)\Phi =0\\[5pt]
\hline 
 \end{array}
$$
}

Properties of the family of quantum algebras
 $U_\zzs(so_{\mma,\mmc}(2,2))$ are now determined by their primitive
generator $P$, since the product $\zzs\pp$ implies that the deformation
parameter $\zzs$ has  dimensions of  {\em length}; hence  we say that
these are  space-type quantum   algebras. Therefore this second quantum
deformation leads to a space discretization of the equation
 (\ref{bf}). Explicitly, if we introduce  the following
differential-difference realization of $U_\zzs(so_{\mma,\mmc}(2,2))$ 
\bea
&&\pp=\partial_x\qquad \hh=\partial_t\nonumber\\[2pt]
 &&\kk= -\mmc t\left(\frac{ {\rm e}^{\zzs\partial_x}-1 }{\zzs}\right)
-\mma x {\rm e}^{-\zzs\partial_x}\partial_t\qquad
 \dd = - x \left(\frac{1-{\rm
e}^{-\zzs\partial_x}}{\zzs}\right)  - t \partial_t  \cr 
&&\cuno=(\mma x^2 {\rm e}^{-2 \zzs\partial_x}+\mmc t^2)\partial_t
 + 2 \mmc x t \left(\frac{1- {\rm e}^{-\zzs\partial_x} }{\zzs}\right)
-\zzs\mma  x {\rm e}^{-2 \zzs\partial_x}\partial_t\cr
&&\cdos=-(\mma x^2 {\rm e}^{-\zzs\partial_x}+\mmc t^2)\left(\frac{ {\rm
 e}^{\zzs\partial_x}-1}{\zzs}\right) - 2 \mma x t\partial_t -\zzs \mma ( t
\partial_t  +  t^2 \partial^2_{t}) 
\label{ed}
\eea
in the Casimir of the deformed subalgebra $\{K,H,P\}$ given by
\be
E_\zzs=\mmc \left( \frac{{\rm e}^{\zzs\pp}-1}{\zzs}\right)^2 -\mma \hh^2  
\label{ee}
\ee
then we obtain a   discretization of the equation
(\ref{bf}) on a uniform space lattice:
\be
E_\zzs\Phi(x,t)=0\quad
\Longrightarrow\quad \left\{\mmc \left(
\frac{{\rm e}^{\zzs\partial_x}-1}{\zzs}\right)^2 -
\mma \frac{\partial^2}{\partial t^2}
\right\}\Phi(x,t)=0 .
\label{ef}
\ee
 The quantum algebra $U_\zzs(so_{\mma,\mmc}(2,2))$ is the symmetry
algebra of this equation as the operators (\ref{ed}) satisfy
\bea
&& [E_\zzs,X]=0\quad \mbox{for}\quad X\in\{\kk,\pp,\hh\}
\qquad
[E_\zzs,\dd]=-2E_\zzs \cr
&&[E_\zzs,\cuno]= 4 \mmc t E_\zzs\qquad
[E_\zzs,\cdos]= -4\mma (x +\zzs + \zzs t \partial_t) E_\zzs .
 \label{eg}
\eea

 To unfold the relationship between these discrete space symmetries and
the results obtained in \cite{Javier} from a Lie symmetry approach  we
consider the twist map for $U_\zzs(so_{\mma,\mmc}(2,2))$ defined by
\bea
&&{\cal P}=\frac{{\rm e}^{\zzs\pp}-1}{\zzs}\qquad
{\cal H}=\hh\qquad {\cal K}=\kk \qquad {\cal D}=\dd \cr
&&{\cal C}_1=\cuno\qquad
{\cal C}_2=\cdos + \zzs \mma \dd^2
\label{eh}
\eea
 that gives rise to the classical commutators (\ref{ba}) with the
coproduct given by 
\be
\begin{array}{l}
\displaystyle{ \Delta({\cal P})=1\otimes {\cal P} + {\cal P} \otimes 1 +
\zzs {\cal P}\otimes {\cal P}}\qquad
\displaystyle{\Delta({\cal H})=1\otimes {\cal H} + {\cal H} \otimes 1 +
\zzs {\cal H}\otimes {\cal P}}\\[5pt]
\displaystyle{\Delta({\cal D})=1\otimes {\cal D} + {\cal D}\otimes 
\frac{1}{1+\zzs{\cal P}}}\qquad
\displaystyle{\Delta({\cal K})=1\otimes {\cal K} + {\cal K} \otimes 1 -
\zzs\mma {\cal D}\otimes \frac{{\cal H}}{1+ \zzs{\cal P}}}\\[8pt]
 \displaystyle{\Delta({\cal C}_1)=1\otimes {\cal C}_1 + {\cal C}_1
\otimes 
\frac{1}{1+\zzs{\cal P}}- 2 \zzs {\cal D} \otimes
 \frac{1}{1+\zzs{\cal P}}\,{\cal K}  +\zzs^2 \mma ({\cal D}^2+{\cal
D})\otimes
\frac{{\cal H}}{(1+\zzs{\cal P})^2} }\\[8pt]
 \displaystyle{ \Delta({\cal C}_2)=1\otimes {\cal C}_2 + {\cal C}_2
\otimes 
\frac{1}{1+\zzs{\cal P}}+ 2 \zzs \mma {\cal D}\otimes
 \frac{1}{1+\zzs{\cal P}}\,{\cal D}  -\zzs \mma  ({\cal D}^2+{\cal
D})\otimes
\frac{\zzs{\cal P}}{(1+\zzs{\cal P})^2} }.
\end{array}
\label{ei}
\ee
Under the map (\ref{eh}), the realization (\ref{ed})  is transformed into
\be
\begin{array}{l}
{\cal P}=\Delta_x\qquad
{\cal H}=\partial_t\\[5pt]
 {\cal K}=  - \mmc t\Delta_x - \mma x T_x^{-1} \partial_t\qquad
{\cal D} = - x T_x^{-1} \Delta_x  - t \partial_t \\[5pt]
{\cal C}_1 =(\mma x^2 T_x^{-2} +\mmc t^2)\partial_t +
2\mmc  x t T_x^{-1} \Delta_x -\zzs\mma x
T_x^{-2}\partial_t \\[5pt]
 {\cal C}_2=-(\mma x^2 T_x^{-2} +\mmc t^2)\Delta_x -
2 \mma x t T_x^{-1} \partial_t +\zzs \mma x  T_x^{-2} \Delta_x  
\end{array}
\label{em}
\ee
 where  $T_x= {\rm e}^{\zzs\partial_x}$ and  $\Delta_x=({T_x-1
})/{\zzs}$. The  element $E_\zzs$ becomes the undeformed  ${\cal E}$
(\ref{ddgg}), so that the  associated differential-difference equation
keeps the form of (\ref{ef}):
 \be 
(\mmc \Delta^2_x -\mma \partial^2_t)\Phi(x,t)=0 . 
\label{en}
\ee
The  operators (\ref{em}) are symmetries of this equation since 
 $\{{\cal K},{\cal H},{\cal P}\}$ commute  with ${\cal E}$ and   the
remaining ones fulfil
\be
[{\cal E},{\cal D}]=-2{\cal E}\qquad
 [{\cal E},{\cal C}_1]= 4\mmc t {\cal E}\qquad [{\cal E},{\cal C}_2]= -
4\mma x T_x^{-1}{\cal E} .
 \label{eo}
\ee
 These last results  reproduce  those found in  \cite{Javier}  once we
introduce the notation (\ref{dj}) and apply the limits  $k\to \infty$,
$\tau\to 0$ (with   $k\tau =t$) in the symmetries and equation of
\cite{Javier} (that  is,  $\Delta_k\to \partial_t$, $T_k\to 1$,
$\Delta_n=\Delta_x$, $T_n=T_x$).   This in turn means that  
$U_\zzs(so_{s^2,1}(2,2))$ is the quantum   algebra  of symmetries of the
equation (\ref{en}) on a uniform space lattice  with   the deformation
parameter $\zzs$  identified with  
 the space lattice constant and   $t$ as a
continuous variable.  The particular equation (\ref{en}) arising for
each  pair $(\mma,\mmc)$    is written in the table 2.


\sect{`Duality' and higher dimensions}

 At a classical level, a remarkable equivalence between the Lie algebras
in the  family  $so_{\mma,\mmc}(2,2)$ (\ref{ba}) is provided by  the map 
defined by
\be
H\to P\quad P\to H\quad K\to K\quad D\to D\quad C_1\to -C_2\quad
C_2\to -C_1
\label{fa}
\ee
that interchanges the role of the generators $H\leftrightarrow P$ and
$C_1\leftrightarrow C_2$,   thus relating the set of
graded contractions as 
\be
so_{\mma,\mmc}(2,2)\leftrightarrow  so_{\mmc,\mma}(2,2).
\label{fb}
\ee
If the  interchange of the two
 coordinates $x \leftrightarrow t$ is added (so $\partial_x
\leftrightarrow \partial_t$), then this  algebraic equivalence also works 
for the vector field realization (\ref{bd}) and equation (\ref{bf}). This
means that if we consider the classical  Lie algebras and associated
differential equations arranged as in table 1  by applying   the
classical limit $\zzt\to 0$ (also as in table 2 for
$\zzs\to 0$),   this  kind of {\em duality} corresponds to the reflection in
the main diagonal. Thus $so(2,2)$ and $i'iso(1,1)$ have self-dual structures,
 meanwhile for the four Lie algebras $iso(2,1)$ this duality  
interchanges the Weyl subalgebras ${\cal WC}\leftrightarrow{\cal WG}$
(isomorphic at this dimension) and the differential equations
$\partial_t^2\Phi=0 \leftrightarrow
\partial_x^2\Phi=0$ according to the transformation of their 
contraction parameters $(\pm,0)\leftrightarrow (0,\pm)$.

 When either the time- or space-type  quantum deformation is introduced in
 the family $so_{\mma,\mmc}(2,2)$, it can be checked that the map
(\ref{fa}) does not lead to a duality as (\ref{fb})  for a {\em single}
family of quantum algebras. To implement this duality at a quantum
algebra level requires to consider  {\em both} families 
 $U_\zzt(so_{\mma,\mmc}(2,2))$ and $U_\zzs(so_{\mma,\mmc}(2,2))$
simultaneously; then the map (\ref{fa})    can be extended by simply
interchanging both deformation parameters  $\zzt\leftrightarrow \zzs$ in
such a manner that both families of quantum algebras are related as
follows
\be
U_\zzt(so_{\mma,\mmc}(2,2)) \leftrightarrow U_\zzs(so_{\mmc,\mma}(2,2)) .
\label{fc}
\ee
Thus the results presented in the table 1 are transformed into those
given in the table 2, and conversely. For instance, the quantum duality
(\ref{fc}) interchanges the  quantum Weyl subalgebras $U_{\zzt}({\cal
WM})\leftrightarrow U_{\zzs}({\cal WM})$,
$U_{\zzt}({\cal WE})\leftrightarrow U_{\zzs}({\cal WE})$,
$U_{\zzt}({\cal WG})\leftrightarrow U_{\zzs}({\cal WC})$ and
$U_{\zzt}({\cal WC})\leftrightarrow U_{\zzs}({\cal WG})$,
 as well as the  derivatives  $\Delta_t \leftrightarrow \Delta_x$ and
$\partial_t \leftrightarrow \partial_x$. Therefore a byproduct of
(\ref{fc}) is that the expressions (\ref{ck}) and (\ref{cm}) become a
matrix realization and an $R$-matrix for $U_\zzs(so_{\mma,\mmc}(2,2))$ 
once  the map (\ref{fa}) has been applied   together with the
replacements $\mma\leftrightarrow \mmc$ and
$\zzt\leftrightarrow \zzs$.

Consequently, both families   of quantum algebras are algebraically
equivalent at this $(1+1)$ dimension. In spite of this fact, we consider
that  the explicit results concerning both families are necessary not
only because from a physical viewpoint they have a different
interpretation and allow us to exhibit the duality clearly, but also
because they  indicate the way to rise to higher dimensions. In this
sense, the $(1+1)$D case is somehow exceptional due to the symmetric role
that   the generators $H$ and $P$  (respectively, the coordinates $t$ and
$x$) play.

We expect that a similar procedure to the one presented in this paper
would enable to construct  quantum deformations for the next dimensions
(keeping the classical $r$-matrices (\ref{ca}) and (\ref{ea}) as the
seeds of the deformations), particularly for the $(3+1)$D case. The
possible quantum  $so(4,2)$ algebras generalizing $U_\zzt(so(2,2))$ and
$U_\zzs(so(2,2))$ would be interpreted as quantum deformations of the
conformal algebra of the  $(3+1)$D Minkowkskian spacetime giving rise to
discretizations of the  $(3+1)$D wave equation as
\be
\begin{array}{ll}
 U_\zzt(so(4,2)):&\quad
(\partial^2_x+\partial^2_y+\partial^2_z-\Delta^2_t )\Phi =0\cr
U_\zzs(so(4,2)):&\quad
(\Delta^2_x+\partial^2_y+\partial^2_z-\partial^2_t)\Phi =0
\end{array}
\label{fd}
\ee
and   fulfilling a sequence of Hopf subalgebras embeddings such as
\be
U_\zzt(sl(2,\R)) \simeq U_\zzt(so(2,1))\subset U_\zzt(so(2,2))\subset
U_\zzt(so(3,2))\subset U_\zzt(so(4,2))\dots
\label{fe}
\ee
We have achieved here the first embedding. In this context we remark
that in \cite{Kulish} (see also  references therein)   a systematic
construction of a chain of twists applied to the universal envelopings of
the semisimple Lie algebras leading to sequences similar to (\ref{fe})
has been introduced.  Furthermore the  structures (\ref{fd}) would be the
cornerstone of a scheme  of contractions   leading to different quantum
deformations of the  algebras  $so(5,1)$, $so(3,3)$, $iso(4,1)$, 
$iso(3,2)$,\dots as well as of their associated differential-difference
equations.

To end with we wish to point out that the quantum $iso(2,1)$ algebras we
have obtained can also be interpreted in a kinematical framework as
quantum  deformations of the $(2+1)$D Poincar\'e algebra by using a {\em 
null-plane  basis} \cite{LS} with
generators $\{P_+,P_1,P_-,E_1,F_1,K_2\}$. If we take, for instance, the
Poincar\'e algebra with contraction parameters $(\mma,\mmc)=(0,+1)$, then
the relationship between the null-plane generators and the conformal ones
is given by
\be
\begin{array}{lll}
P_+=\frac {1}{\sqrt{2}}\, P&\quad P_1=K&\quad P_-=-\frac {1}{\sqrt{2}}\,
C_2\cr
 E_1=-\frac {1}{\sqrt{2}}\, H&\quad F_1=\frac {1}{\sqrt{2}}\, C_1
&\quad K_2 =D .
\end{array}
\ee
This change of basis gives rise to  two inequivalent quantum Poincar\'e
algebras: $U_\zzt(iso(2,1))\supset U_\zzt(so(2,1))$ with $E_1$ primitive
and $U_\zzs(iso(2,1))\supset U_\zzs(iso(1,1))$ with $P_+$   primitive.
 These non-standard deformations are different from the so called
null-plane quantum Poincar\'e algebra \cite{nulla,nullb}, the  underlying
classical $r$-matrix  of which reads  $r=2z (K_2\wedge P_+ + E_1\wedge
P_1)$ in the $(2+1)$D case.


\section*{Acknowledgment}
 This work was partially supported  by  Junta de Castilla y
Le\'on, Spain  (Project   CO2/399).


\end{document}